\documentclass[10pt,a4paper,final]{amsart}

\usepackage{amssymb}
\def\intt{\operatorname{int}}
\def\CC{\mathbb C}
\def\MT{\mathbb T}
\def\MX{\mathbb X}
\def\NN{\mathbb N}
\def\QQ{\mathbb Q}
\def\CM{\mathcal M}
\def\CO{\mathcal O}
\def\halfskip{\vskip 10pt plus 1pt minus 1pt}
\let\phi=\varphi
\def\PSH{\mathcal P\mathcal S\mathcal H}

\def\tdots{\!\times\!\dots\!\times\!}
\def\theequation{\arabic{chapter}.\arabic{section}.\arabic{equation}}
\def\too{\longrightarrow}
\def\wdht{\widehat}
\def\wdtl{\widetilde}
\def\Delta{\varDelta}
\def\Gamma{\varGamma}
\def\Omega{\varOmega}
\catcode`@=11
\def\@opargbegintheorem#1#2#3{%
\item[\hskip\labelsep \theorem@headerfont #1\ #2]\rom{(#3)}.}
\def\th@mytheorem{%
  \let\thm@indent\noindent 
  \thm@headfont{\bfseries}
  \itshape 
}
\def\th@myremark{%
  \let\thm@indent\noindent 
  \thm@headfont{\bfseries}
}
\catcode`@=12

\theoremstyle{mytheorem}
\newtheorem{Theorem}{Theorem}[section]
\newtheorem{Corollary}[Theorem]{Corollary}

\theoremstyle{myremark}
\newtheorem{Definition}[Theorem]{Definition}
\newtheorem{Remark}[Theorem]{Remark}
\catcode`@=11
\tagsleft@false
\def\@eqnnum{\tagform@\theequation}
\def\@makefnmark{\hbox{$\left(^{\@thefnmark}\right)$\;}}
\@mparswitchfalse
\renewcommand\qed{\RIfM@\else\unskip\nobreak\fi\hfill\qedsymbol}
\def\partrunhead#1#2#3{\small%
  \@ifnotempty{#2}{{\ignorespaces#1 #2\unskip}\@ifnotempty{#3}{. }}%
  \def\@tempa{#3}%
  \ifx\@empty\@tempa\else\@tempa\fi}

\def\@cite#1#2{{%
 \m@th\upshape\mdseries[{#1}{\if@tempswa, #2\fi}]}}
\@ifundefined{cite }{%
  \expandafter\let\csname cite \endcsname\cite
  \edef\cite{\@nx\protect\@xp\@nx\csname cite \endcsname}%
}{}
\catcode`@=12

\begin{document}

\title[An extension theorem for separately meromorphic functions]
{An extension theorem for separately meromorphic functions with pluripolar
singularities}

\author{Marek Jarnicki}
\address{Jagiellonian University, Institute of Mathematics\newline
\indent Reymonta 4, 30-059 Krak\'ow, Poland}
\email{jarnicki@im.uj.edu.pl}
\thanks{The first author was supported in part by the KBN grant
No.~5 P03A 033 21. The paper was finished during the RiP stay of both authors
at Mathematisches Forschungsinstitut in Oberwolfach. We like to thank these
institutions.}

\author{Peter Pflug}
\address{Carl von Ossietzky Universit\"at Oldenburg, Fachbereich Mathematik\newline
\indent Postfach 2503, D-26111 Oldenburg, Germany}
\email{pflug@{mathematik.uni-oldenburg.de}}



\begin{abstract}
Let $D_j\subset\CC^{n_j}$ be a pseudoconvex domain and let $A_j\subset D_j$
be a locally pluriregular set, $j=1,\dots,N$. Put
$$
X:=\bigcup_{j=1}^N A_1\tdots A_{j-1}\times D_j\times A_{j+1}\tdots A_N.
$$
Let $M\subset X$ be relatively closed.
For any $j\in\{1,\dots,N\}$ let
$\Sigma_j$ be the set of all $(z',z'')\in(A_1\tdots A_{j-1})
\times(A_{j+1}\tdots A_N)$ such that the fiber
$M_{(z',\cdot,z'')}:=\{z_j\in\CC^{n_j}: (z',z_j,z'')\in M\}$ is not pluripolar.
Assume that $\Sigma_1,\dots,\Sigma_N$ are pluripolar. Put
\begin{multline*}
X':=\bigcup_{j=1}^N\{(z',z_j,z'')\in(A_1\tdots A_{j-1})\times D_j
\times(A_{j+1}\tdots A_N):\\
(z',z'')\notin\Sigma_j\}.
\end{multline*}
Then (Theorem \ref{OMT}) there exists a relatively closed pluripolar
subset $\wdht M\subset\wdht X$ of the `envelope of holomorphy'
$\wdht X$ of $X$ such that:

$\bullet$ $\wdht M\cap X'\subset M$,

$\bullet$ every function $f$ separately meromorphic on $X\setminus M$
(Definition \ref{SMF}) extends to a (uniquely determined) function
$\wdht f$ meromorphic on $\wdht X\setminus\wdht M$,

$\bullet$ if $f$ is separately holomorphic on $X\setminus M$, then $\wdht f$
is holomorphic on $\wdht X\setminus\wdht M$, and

$\bullet$ $\wdht M$ is singular with respect to the family of all functions
$\wdht f$.

\noindent The case of separately holomorphic functions was solved in
\cite{JarPfl2002b}. In the case where $N=2$, $M=\varnothing$, the above result will
be strengthened in Theorem \ref{MNT}.

\end{abstract}


\maketitle

\section{Introduction. Main results.}

\halfskip

We keep the main notation from \cite{JarPfl2002b}:

$\bullet$ Let $N\in\NN$, $N\geq2$, and let
$\varnothing\neq A_j\subset D_j\subset\CC^{n_j}$,
where $D_j$ is a domain, $j=1,\dots,N$. We define an {\it $N$--fold cross}
\begin{align*}
X&=\MX(A_1,\dots,A_N;D_1,\dots,D_N)\\
&:=\bigcup_{j=1}^NA_1\tdots A_{j-1}\times D_j\times A_{j+1}\tdots A_N
\subset\CC^{n_1+\dots+n_N}=\CC^n.
\end{align*}

$\bullet$ For an open set $\Omega\subset\CC^n$ and $A\subset\Omega$ let
$$
h_{A,\Omega}:=\sup\{u:\; u\in\PSH(\Omega),\;u\leq1 \text{ on }
\Omega,\; u\leq0 \text{ on } A\},
$$
where $\PSH(\Omega)$ is the set of all functions plurisubharmonic on
$\Omega$. Put
$$
\omega_{A,\Omega}:=\lim_{k\to+\infty}h^\ast_{A\cap\Omega_k,\Omega_k},
$$
where $(\Omega_k)_{k=1}^\infty$ is a sequence of relatively compact open sets
$\Omega_k\subset\Omega_{k+1}\Subset\Omega$ with
$\bigcup_{k=1}^\infty\Omega_k=\Omega$ ($h^\ast$ denotes the upper
semicontinuous regularization of $h$).

$\bullet$ For an $N$--fold cross $X=\MX(A_1,\dots,A_N;D_1,\dots,D_N)$ put
$$
\wdht X:=\{(z_1,\dots,z_N)\in D_1\tdots D_N:
\sum_{j=1}^N\omega_{A_j,D_j}(z_j)<1\}.
$$

$\bullet$ We say that a subset $\varnothing\neq A\subset\CC^n$ is {\it locally
pluriregular} if $h^\ast_{A\cap\Omega,\Omega}(a)=0$ for any $a\in A$ and for
any open neighborhood $\Omega$ of $a$.

$\bullet$ Suppose that $S_j\subset (A_1\tdots A_{j-1})
\times(A_{j+1}\tdots A_N)$, $j=1,\dots,N$.
Define the {\it generalized $N$--fold cross}
\begin{multline*}
T=\MT(A_1,\dots,A_N;D_1,\dots,D_N;S_1,\dots,S_N)
:=\bigcup_{j=1}^N\{(z',z_j,z'')\\
\in (A_1\tdots A_{j-1})\times D_j\times(A_{j+1}\tdots A_N): (z',z'')\notin S_j\}.
\quad\footnotemark
\end{multline*}

\footnotetext{Observe that $\MX(A_1,\dots,A_N;D_1,\dots,D_N)=
\MT(A_1,\dots,A_N;D_1,\dots,D_N;\varnothing,\dots,\varnothing)$. Moreover,
if $N=2$, then $\MT(A_1,A_2;D_1,D_2;S_1,S_2)=
\MX(A_1\setminus S_2,A_2\setminus S_1;D_1,D_2)$.}

$\bullet$
Let $M\subset T$ be a relatively
closed set. We say that a function $f:T\setminus M\too\CC$ is {\it separately
holomorphic} ($f\in\CO_s(T\setminus M)$) if for any $j\in\{1,\dots,N\}$ and
$(a',a'')\in(A_1\tdots A_{j-1})\times(A_{j+1}\tdots A_N)\setminus S_j$ the function
$f(a',\cdot,a'')$ is holomorphic in the open set
$D_j\setminus M_{(a',\cdot,a'')}$, where
$M_{(a',\cdot,a'')}:=\{z_j\in\CC^{n_j}: (a',z_j,a'')\in M\}$ \footnote{Observe
that the above condition is empty if $M_{(a',\cdot,a'')}=D_j$.}\!. Notice that
the definition applies to the case where $T=X$ is an $N$--fold cross
($S_1=\dots=S_N=\varnothing$).

\halfskip

The following general extension theorem for separately holomorphic functions
with singularities was proved in \cite{JarPfl2002a} and \cite{JarPfl2002b}.

\begin{Theorem}\label{MT}
Let $D_j\subset\CC^{n_j}$ be a pseudoconvex domain,
let $A_j\subset D_j$ be a locally pluriregular set, $j=1,\dots,N$, and let
$M\subset X$ be a relatively closed subset of the $N$--fold cross
$X:=\MX(A_1,\dots,A_N;D_1,\dots,D_N)$.
Assume that for each $j\in\{1,\dots,N\}$ the set
$\Sigma_j=\Sigma_j(M)$ of all points
$(z',z'')\in(A_1\tdots A_{j-1})\times(A_{j+1}\tdots A_N)$ such that the fiber
$M_{(z',\cdot,z'')}$ is not pluripolar is pluripolar. Put
$$
X'=X'(M):=\MT(A_1,\dots,A_N;D_1,\dots,D_N;\Sigma_1,\dots,\Sigma_N).
$$
Then there exists a relatively closed pluripolar
set $\wdht M\subset\wdht X$ such that:

$\bullet$ $\wdht M\cap X'\subset M$,

$\bullet$ for every $f\in\CO_s(X\setminus M)$ there exists exactly one
$\wdht f\in\CO(\wdht X\setminus\wdht M)$ with $\wdht f=f$ on
$X'\setminus M$,

$\bullet$ $\wdht M$ is singular with respect to the family
$\{\wdht f: f\in\CO_s(X\setminus M)\}$.

In particular, $\wdht X\setminus\wdht M$ is the envelope of holomorphy of
$X\setminus M$ with respect to the space of separately holomorphic functions.

Moreover:

{\rm (a)} if $M$ is pluripolar,
then $\Sigma_1,\dots,\Sigma_N$ are pluripolar \footnote{And, consequently, the
assumption of the theorem is always satisfied for pluripolar sets.}\!,

{\rm (b)} if $M=X\cap\wdtl M$, where $\wdtl M$ is an analytic
subset of an open connected neighborhood of $X$,
then $\wdht M$ is analytic,

{\rm (c)} if $M=X\cap\wdtl M$, where $\wdtl M$ is an
analytic subset of $\wdht X$, then $\wdht M$ is the union of all pure
$(n-1)$--dimensional irreducible components of $\wdtl M$
\footnote{In particular,
$M=\varnothing \Longrightarrow \wdht M=\varnothing$.}\!.
\end{Theorem}

Some special cases of the above theorem were studied by many authors ---
see the references in \cite{JarPfl2002b}.

It is known that the envelope of holomorphy (of any Riemann domain over
$\CC^n$) coincides with the envelope of meromorphy (cf.~\cite{JarPfl2000},
Th.~3.6.6). Thus it is natural to conjecture that in the above situation
the domain $\wdht X\setminus\wdht M$ is also the envelope of meromorphy of
$X\setminus M$ with respect to separate meromorphic functions. The case
$M=\varnothing$ was studied in \cite{Sak1957}, \cite{Kaz1976}, \cite{Kaz1978},
\cite{Kaz1984}, \cite{Shi1986}, and \cite{Shi1989}.

\begin{Definition}\label{SMF}
Let $T=\MT(A_1,\dots,A_N;D_1,\dots,D_N;S_1,\dots,S_N)$ be a generalized
$N$--fold cross. Let $M\subset T$, $S\subset T\setminus M$
be relatively closed. We say that a function
$f:(T\setminus M)\setminus S\too\CC$ is {\it separately meromorphic on
$T\setminus M$} ($f\in\CM_s(T\setminus M)$)
if for any $j\in\{1,\dots,N\}$ and $(a',a'')\in (A_1\tdots A_{j-1})
\times(A_{j+1}\tdots A_N)\setminus S_j$ with
$(M\cup S)_{(a',\cdot,a'')}\neq D_j$, there exists a function $\wdtl{f(a',\cdot,a'')}\in\CM(D_j\setminus M_{(a',\cdot,a'')})$
such that $\wdtl{f(a',\cdot,a'')}=f(a',\cdot,a'')$ on $D_j\setminus(M\cup S)
_{(a',\cdot,a'')}$.
\end{Definition}
Observe that $f\in\CO_s(T\setminus(M\cup S))$
\footnote{Note that $M\cup S$ is relatively closed in $T$.}\!.

\halfskip

The main results of the paper are the following two theorems.

\begin{Theorem}\label{OMT}
Let $(A_j,D_j)_{j=1}^N$, $X$, $M$, and $\wdht M$ be as in Theorem \ref{MT}.
Let $S\subset X\setminus M$ be relatively closed and
let $f:(X\setminus M)\setminus S\too\CC$ be a separate meromorphic function
on $X\setminus M$ such that

{\rm (*)} the sets $\Sigma_1(S),\dots,\Sigma_N(S)$ are pluripolar.

Put $Q_f:=M\cup S$. Then there exists exactly one
$\wdht f\in\CM(\wdht X\setminus\wdht M)$ such that:

$\bullet$ $\wdht f\in\CO(\wdht X\setminus\wdht Q_f)$, where
the set $\wdht Q_f$ is constructed via Theorem \ref{MT} (in the same way
as $\wdht M$ for $M$) \footnote{Note that $\wdht M\subset\wdht Q_f$.}\!,

$\bullet$ $\wdht f=f$ on $X'_f\setminus Q_f$, where
$$
X'_f:=\MT(A_1,\dots,A_N;D_1,\dots,D_N;\Sigma_1(Q_f),\dots, \Sigma_N(Q_f)).
$$

Consequently, the envelope of $X\setminus M$ with respect to separately
meromorphic functions satisfying {\rm (*)} coincides with its
envelope of separate holomorphy.
\end{Theorem}

In the case where $N=2$, $M=\varnothing$, the above result may be
strengthened as follows.

\begin{Theorem}\label{MNT}
Let $D\subset\CC^p$, $G\subset\CC^q$ be pseudoconvex domains, let
$\varnothing\neq A\subset D$, $\varnothing\neq B\subset G$ be locally
pluriregular sets, and let
$$
X:=\MX(A,B;D,G)=(A\times G)\cup(D\times B).
$$
Let $S\subset X$ be a relatively closed set. Assume that:

{\rm (1.4.1)} for every $(a,b)\in A\times B$ we have
$\intt_{\CC^q}S_{(a,\cdot)}=\varnothing$,
$\intt_{\CC^p}S_{(\cdot,b)}=\varnothing$,

{\rm (1.4.2)} $A\times B\subset\overline{(A\times B)\setminus S}$
\footnote{In particular, for every
$(a,b)\in A\times B$ and for every neighborhood $U\subset\CC^p\times\CC^q$
of $(a,b)$ the set $(A\times B)\cap U\setminus S$ is not pluripolar.}\!,

there exist exhaustions $(D_j)_{j=1}^\infty$ and
$(G_j)_{j=1}^\infty$ of $D$ and $G$, respectively, such that:

{\rm (1.4.3)} $D_j$, $G_j$ are relatively closed pseudoconvex subdomains of
$D$ and $G$, respectively,

{\rm (1.4.4)} $A_j:=A\cap D_j\neq\varnothing$, $B_j:=B\cap G_j\neq\varnothing$,

{\rm (1.4.5)} for every $(a,b)\in A_j\times B_j$ we have $B_j\setminus S_{(a,\cdot)}\neq\varnothing$,
$A_j\setminus S_{(\cdot,b)}\neq\varnothing$,
$j=1,2,\dots$.

Then for every function $f:X\setminus S\too\CC$ which is separately
meromorphic on $X$ there exists a function $\wdht f\in\CM(\wdht X)$ such
that $\wdht f=f$ on $X\setminus S$.
\end{Theorem}


\halfskip

\section{Auxiliary results.}

\halfskip

\begin{Remark}\label{ArmGar}
(a) (\cite{Kli1991}, Corollary 4.8.4) If $A, B\subset\CC^n$ are plurithin
at a point $a\in\CC^n$ \footnote{We say that a set $A\subset\CC^n$ is
{\it plurithin} at a point $a\in\CC^n$ if either $a\notin\overline A$ or
$a\in\overline A$ and $\limsup_{A\setminus\{a\}\ni z\to a}u(z)<u(a)$
for a function $u$ plurisubharmonic in a neighborhood of $a$.}\!,
then $A\cup B$ is plurithin at $a$.

(b) (\cite{ArmGar2001}, Th. 7.2.2) Every polar set $P\subset\CC$ is
thin at any point $a\in\CC$.

(c) If $A\subset\CC$ is not thin at a point $a\in\overline A$, then
for any polar set $P\subset\CC$, the set $A\setminus P$ is not thin
at $a$ ((c) follows directly from (a) and (b)).

(d) If $A\subset\CC^n$ is locally pluriregular at a point $a\in\overline A$,
then $A$ is not plurithin at $a$.
If $A\subset\CC$ is not thin at a point $a\in\overline A$, then  $A$ is
locally regular at $a$.

Indeed,
suppose that $A\subset\CC^n$ is locally pluriregular at $a$ and
$$
\limsup_{A\setminus\{a\}\ni z\to a}u(z)<c<u(a)
$$
for some $u\in\PSH(V)$, where $V$ is an open neighborhood of $a$. We may assume that
$u\leq0$ on $V$. Take an open neighborhood $U\subset V$ of $a$
such that $u<c$ on $(A\setminus\{a\})\cap U$.
Put $v:=\frac{u}{-c}+1$.
Then $v\leq1$ on $U$ and $v\leq0$ on $(A\setminus\{a\})\cap U$. Hence
$v\leq h^\ast_{(A\setminus\{a\})\cap U,U}=h^\ast_{A\cap U,U}$ on $U$.
In particular, $0=v(a)=\frac{u(a)}{-c}+1<0$; contradiction.

Now, suppose that $A\subset\CC$ is not thin at $a$ and $h^\ast_{A\cap U,U}(a)>0$
for some neighborhood $U$ of $a$. Let $P\subset U$ be a polar set such
that $h^\ast_{A\cap U,U}=h_{A\cap U,U}$ on $U\setminus P$ (cf.~\cite{JarPfl2000}
Th.~2.1.41).
In particular, $h^\ast_{A\cap U,U}=0$ on $A\setminus P$.
By (c), the set $A\setminus P$ is not thin at $a$. Hence
$0<h^\ast_{A\cap U,U}(a)=\limsup_{A\setminus P\ni z\to a}
h^\ast_{A\cap U,U}(z)=0$; contradiction.

(e) (\cite{ArmGar2001}, Th.~7.3.9) If $A\subset\CC$ is thin
at a point $a\in\overline A$, then there is a sequence $r_k\searrow0$ such that $\{z\in A:|z-a|=r_k\}=
\varnothing$, $k=1,2,\dots$.

(f) (\cite{BedTay1982}, Corollary 10.5) For a non-pluripolar set $A\subset\CC^n$
let $A^\ast$ denote the set
of all $a\in\overline A$ such that $A$ is locally pluriregular at $a$. Then
$A\setminus A^\ast$ is pluripolar.
\end{Remark}


\halfskip

\section{Corollaries from Theorem 1.4.}

\halfskip

Let $E$ denote the unit disc. For $a\in\CC^k$, $r>0$, let
$\Delta_a(r)=\Delta_a^k(r)$ be the polydisc with center at $a$ and the
radius $r$.

\begin{Corollary}[Cf. \cite{Sak1957}]\label{Sak}
Let $S\subset E\times E$ be a relatively closed set such that:

$\bullet$ $\intt S=\varnothing$,

$\bullet$ for every domain  $U\subset E\times E$ the set $U\setminus S$ is
connected \footnote{We shortly say that $S$ does not separate domains.}\!.

Let $A$ (resp. $B$) denote the set of all $a\in E$ (resp. $b\in E$) such that
$\intt_{\CC}S_{(a,\cdot)}=\varnothing$ (resp. $\intt_{\CC}S_{(\cdot,b)}=
\varnothing$). Put $X:=\MX(A,B;E,E)=(A\times E)\cup(E\times B)$.

Then for every function $f:X\setminus S\too\CC$ which is separately
meromorphic on $X$, there exists an $\wdht f\in\CM(E\times E)$ such that
$\wdht f=f$ on $X\setminus S$.
\end{Corollary}

\begin{Remark}
Notice that the original proof of the above result is not correct:
the proof of Theorem 1 in \cite{Sak1957} contains an essential gap.
Namely, on p.~78 the author claims that for any domain $U\subset E\times E$
the set $(A\times B)\cap U\setminus S$ contains an open polydisc.
The following example shows that this is in general impossible.

Let $(\QQ+i\QQ)\cap E=\{q_1, q_2,\dots\}$, $S:=\bigcup_{k=1}^\infty
\{q_k\}\times\overline\Delta_{1-1/k}(1/k^2)\subset E\times E$.
Then $S$ satisfies all the assumptions of Corollary \ref{Sak} but in this case
the interior of $A=E\setminus(\QQ+i\QQ)$ is empty.
\end{Remark}

\begin{proof} First we check that the sets $A$ and $B$ are not thin at
any point of $E$ (in particular, they are dense in $E$).

Indeed, suppose that $A$ is thin at a point $a\in E$. By Remark
\ref{ArmGar}(e), there exist a circle $C\subset E$ such that
$C\cap A=\varnothing$. Using a Baire category argument, we conclude that
there exist a non-empty open arc $\Gamma\subset C$ and an open disc
$\Delta\subset E$ such that the $3$--dimensional real surface
$\Gamma\times\Delta$ is contained in $S$. Hence, since $S$ is nowhere dense
and does not separate domains, we get a contradiction.

Consequently, by Remark \ref{ArmGar}(d), the sets $A$ and $B$ are locally
regular and $h^\ast_{A,E}=h^\ast_{B,E}=0$. In particular, $\wdht X=E\times E$.

Now, using the fact that $A$ and $B$ are dense in $E$, one can easily check
that all the assumptions of Theorem \ref{MNT} ($D=G=E$) are satisfied
with arbitrary exhaustions $D_j:=\Delta_0(r_j)$, $G_j:=\Delta_0(r_j)$,
$0<r_j\nearrow1$, which satisfy condition (1.4.3).
\end{proof}

\begin{Remark} (a) E. Sakai claims in \cite{Sak1957} that also the following
$n$--dimensio\-nal version of Corollary \ref{Sak} is true. We do not know how
to prove it.

{\it Let $S\subset E^n$ be relatively closed such that $\intt S=\varnothing$
and $S$ does not separate domains.
Let $f:E^n\setminus S\too\CC$ be such that for any $j\in\{1,\dots,n\}$ and
for any $(a',a'')\in E^{j-1}\times E^{n-j}$ for which
$\intt_{\CC}S_{(a',\cdot,a'')}=\varnothing$,
the function $f_{(a',\cdot,a'')}$ extends meromorphically to $E$
\footnote{That is, $f$ is separately meromorphic on the $n$--fold
generalized cross $T:=\MT(E,\dots,E;E,\dots,E;S_1,\dots,S_n)$, where
$S_j$ denote the set of all $(a',a'')\in E^{j-1}\times E^{n-j}$ for which
$\intt_{\CC}S_{(a',\cdot,a'')}\neq\varnothing$, $j=1,\dots,n$;
cf. Definition \ref{SMF}.}\!.
Then $f$ extends meromorphically to $E^n$.
}

In particular, we would like to ask whether for any set $A\subset E^k$
which is plurithin at $0\in E^k$ there exists a non-empty relatively open
subset $\Gamma$ of a real hypersurface such that
$\Gamma\subset E^k\setminus A$ (cf.~the proof of Corollary \ref{Sak}).

(b) We also do not know whether the following generalization of Corollary
\ref{Sak} is true.

{\it Let $D\subset\CC^p$, $G\subset\CC^q$ be pseudoconvex domains and let
$S\subset D\times G$ be a relatively closed set such that $\intt S=\varnothing$
and $S$ does not separate domains. Let $A$ (resp. $B$) denote the set of all
$a\in D$ (resp. $b\in G$) such that $\intt_{\CC^q}S_{(a,\cdot)}=\varnothing$
(resp. $\intt_{\CC^p}S_{(\cdot,b)}=\varnothing$). Put
$X:=\MX(A,B;D,G)=(A\times G)\cup(D\times B)$.
Then for every function $f:X\setminus S\too\CC$ which is separately
meromorphic on $X$, there exists an $\wdht f\in\CM(D\times G)$ such that
$\wdht f=f$ on $X\setminus S$.
}
\end{Remark}

\begin{Corollary}[Cf. \cite{Shi1989}, Th. 2]\label{Shi}
Let $D,G,A,B,X$ be as in Theorem \ref{MNT}. Assume that $S\subset X$
is a relatively closed set such that

$\bullet$ the set $D\setminus A$ is of zero Lebesgue measure,

$\bullet$ for every $a\in A$ the fiber $S_{(a,\cdot)}$ is pluripolar,

$\bullet$ for every $b\in B$ the fiber $S_{(\cdot,b)}$ is of zero Lebesgue
measure.

Then for every function $f:X\setminus S\too\CC$ which is separately
meromorphic on $X$, there exists an
$\wdht f\in\CM(D\times G)$ such that $\wdht f=f$ on $X\setminus S$.
\end{Corollary}

\begin{proof} One can easily check that all the assumptions of
Theorem \ref{MNT} are satisfied (with arbitrary exhaustions satisfying
(1.4.3--4)). It remains to observe that $h^\ast_{A,D}\equiv0$ (because
$h^\ast_{A,D}=0$ on $A$ and the set $D\setminus A$ is of zero measure).
Hence $\wdht X=D\times G$.
\end{proof}


\halfskip

\section{Rothstein theorem.}

\halfskip

\begin{Theorem}[Cf. \cite{Rot1950}]\label{Rots}
Let $f\in\CM(E^p\times E^q)$. Assume that $A\subset E^p$ be a locally
pluriregular set such that for any $a\in E^p$ we have
$(P_f)_{(a,\cdot)}\neq E^q$, where $P_f$ denote the polar set of $f$,
i.e. $P_f$ is the union of the set of all poles of $f$ and the set of
all indeterminancy points of $f$ \footnote{Note that $P_f$ is analytic and
$f\in\CO(E^p\times E^q\setminus P_f)$.}\!. Let $G\subset\CC^p$ be a domain
such that $E^q\subset G$. Assume that for every $a\in A$ the function
$f(a,\cdot)$ extends meromorphically to $G$. Then there exists an open
neighborhood $\Omega$ of $(E^p\times E^q)\cup(A\times G)$ and a function
$\wdht f\in\CM(\Omega)$ such that $\wdht f=f$ on $E^p\times E^q$.
\end{Theorem}

We present a sketch of the proof.

(1) The case where $A=E^p$
\footnote{Observe that if $A=E^p$, then we
have to prove that $f$ extends meromorphically to $E^p\times G$.}\!,
$q=1$, $G=\Delta_0(R)$ ($R>1$), and $f\in\CO(E^p\times E)$:

The proof may be found for instance in \cite{Siu1974}.

\halfskip

(2) The case where $A=E^p$, $q=1$, and $G=\Delta_0^q(R)$:

Recall that $(P_f)_{(a,\cdot)}\neq E^q$ for any $a\in E^p$, and therefore, for
any $a\in E^p$ there exists a $b\in E^q$ such that $f$ is holomorphic in a
neighborhood of $(a,b)$. By applying locally (1), we get the required result.

\halfskip

(3) The case where $A=E^p$ and $G=\Delta_0^q(R)$:

Let $R_0$ denote the radius of the maximal polydisc $\Delta_0^q(R_0)$ such
that $f$ extends meromorphically to $E^p\times\Delta_0^q(R_0)$. We only need
to show that $R_0\geq R$. Obviously $R_0\geq1$. Suppose that $R_0<R$.

Let $S_q$ be the set of all $(z,w')\in E^p\times\Delta_0^{q-1}(R_0)$ such that
$(P_f)_{(z,w',\cdot)}=E$. It is well known that $S_q$ is an analytic subset of
$E^p\times\Delta_0^{q-1}(R_0)$. Moreover, our assumptions imply that
$S_q\neq E^p\times\Delta_0^{q-1}(R_0)$. Applying locally the Rothstein theorem
to $(E^p\times\Delta_0^{q-1}(R_0)\setminus S_q)\times\Delta_0(R)\subset\CC^{p+q-1}\times\CC$,
we conclude that $f$ extends meromorphically to
$((E^p\times\Delta_0^{q-1}(R_0)\setminus S_q)\times\Delta_0(R))\cup
(E^p\times\Delta_0^q(R_0))$. Observe that, by the
Levi extension theorem (\cite{JarPfl2000}, Prop. 3.4.5),
the envelope of holomorphy of $((E^p\times\Delta_0^{q-1}(R_0)
\times\Delta_0(R))\setminus(S_q\times\Delta_0(R)))\cup
(E^p\times\Delta_0^q(R_0))$ equals $E^p\times\Delta_0^{q-1}(R_0)
\times\Delta_0(R)$. Consequently, the
function $f$ extends meromorphically to
$E^p\times\Delta_0^{q-1}(R_0)\times\Delta_0(R)$.
Repeating the same argument with respect to other variables in $\CC^q$,
we conclude that $f$ extends meromorphically to the domain
$E^p\times H$, where
$$
H=\bigcup_{j=1}^q\Delta_0^{j-1}(R_0)\times\Delta_0(R)\times\Delta_0^{q-j}(R_0).
$$
The envelope of holomorphy of $E^p\times H$
has the form $E^p\times\wdht H$, where $\wdht H$ contains a polydisc
$\Delta_0^q(R'_0)$ with $R'_0>R_0$. Thus $f$ extends meromorphically to
$E^p\times\Delta_0^q(R'_0)$; contradiction --- cf. the proof of Lemma 12
in \cite{JarPfl2002b}.

\halfskip

(4) The case where $A\subset E^p$ is locally pluriregular and
$G=\Delta_0^q(R)$:

For every $z\in E^p$,
let $\rho_f(z)$ denote the radius of the maximal polydisc
$\Delta_0^q(r)$ such that $f(z,\cdot)$ extends meromorphically to
$\Delta_0^q(r)$. Obviously, $\rho_f\geq1$ on $E^p$ and $\rho_f\geq R$ on $A$.

Using (3),  one can easily conclude that $f$ extends
meromorphically to the Hartogs domain
$$
D:=\{(z,w)\in E^p\times\CC^q: |w|<(\rho_f)_\ast(z)\}.
$$
Let $\wdtl f\in\CM(D)$ be the meromorphic extension of $f$.

Moreover, $-\log(\rho_f)_\ast\in\PSH(E^p)$.

Indeed, let $\wdht D$ denote the envelope of holomorphy of $D$. It is known
that $\wdht D\subset E^p\times\CC^q$ is a Hartogs domain with complete
$q$--circled fibers (\cite{JarPfl2000}, Remark 3.1.2(h)). Moreover,
$\wdtl f$ extends meromorphically to $\wdht D$
(\cite{JarPfl2000}, Th. 3.6.6). In particular,
$$
(\rho_f)_\ast(z)=\inf\{\delta_{\wdht D,(0,\xi)}(z,0): \xi\in\CC^q,\; |\xi|=1\},
\quad z\in E^p,
$$
where
$$
\delta_{\wdht D,(0,\xi)}(z,0)=\sup\{r>0: (z,0)+\Delta_0(r)(0,\xi)
\subset\wdht D\}.
$$
Consequently, $-\log(\rho_f)_\ast\in\PSH(E^p)$
(\cite{JarPfl2000}, Th. 2.2.9(iv)).

\halfskip

Thus $-\log(\rho_f)_\ast\in\PSH(E^p)$. Recall that $\rho_f\geq R$ on $A$.
Hence, using the local pluriregularity of $A$, we conclude that
$(\rho_f)_\ast\geq R$ on $A$ \footnote{Suppose that
$h^\ast_{A,E^p}=h_{A,E^p}$ on $E^p\setminus P$, where $P$ is pluripolar.
Put $u:=\frac{-\log(\rho_f)_\ast}{\log R}+1$. Then $u\leq1$ and $u\leq0$ on
$A\setminus P$. Consequently, $u\leq h^\ast_{A\setminus P,E^p}=h^\ast_{A,E^p}$.
In particular, $u\leq0$ on $A$, i.e. $(\rho_f)_\ast\geq R$ on $A$.}\!.
Thus $A\times\Delta_0^q(R)\subset D$,
and therefore $D$ is the required neighborhood.

\halfskip

(5) The general case where $A\subset E^p$ is locally pluriregular and
$G$ is arbitrary:

Fix an $a\in A$.
Let $G_0$ denote the set of all $b\in G$ such that there exist $r_b>0$ and
$f_b\in\CM(\Delta_{(a,b)}(r_b))$, $\Delta_{(a,b)}(r_b)\subset E^p\times G$,
such that:

$\forall_{\alpha\in A\cap\Delta_a(r_b)}: f_b(\alpha,\cdot)=
\wdtl{f(\alpha,\cdot)}$ on $\Delta_b(r_b)$ \footnote{As before,
$\wdtl{f(\alpha,\cdot)}$ denotes the meromorphic extension of
$f(\alpha,\cdot)$.}\!.

Obviously $G_0$ is open, $G_0\neq\varnothing$ ($E^q\subset G_0$). Using
the Rothstein theorem with $G=\Delta_0^q(R)$, one can prove that
$G_0$ is closed in $G$. Thus $G_0=G$.

Moreover, one can also prove that if
$\Delta_{b'}(r_{b'})\cap\Delta_{b''}(r_{b''})\neq\varnothing$, then
$f_{b'}=f_{b''}$ on $\Delta_{(a,b')}(r_{b'})\cap\Delta_{(a,b'')}(r_{b''})$.
This gives a meromorphic extension of $f$ to an open neighborhood of
$\{a\}\times G$. Since $a$ was arbitrary, we get the required neighborhood
$\Omega$.

\halfskip

The proof of the Rothstein theorem is completed.


\halfskip

\section{Proof of Theorem 1.3.}

\halfskip

Fix a function
$f\in\CM_s(X\setminus M)\cap\CO_s(X\setminus Q_f)$.
By Theorem \ref{MT} there exists exactly one
$\wdht f\in\CO(\wdht X\setminus\wdht Q_f)$ with $\wdht f=f$ on
$X'_f\setminus Q_f$.
It remains to prove that $\wdht f\in\CM(\wdht X\setminus\wdht M)$.

It is sufficient to prove that $\wdht f\in\CM(\Omega\setminus\wdht M)$,
where $\Omega\subset\wdht X$ is an open neighborhood of $X'_f$.

Indeed, by virtue of Lemma 9 from \cite{JarPfl2002b} and the Chirka theorem
(cf. \cite{Chi1993}, see also \cite{JarPfl2002b}, Th.~6), the envelope of
holomorphy of $\Omega\setminus\wdht M$ coincides with
$\wdht X\setminus\wdht M$. Consequently,
the function $\wdht f$ extends meromorphically to $\wdht X\setminus\wdht M$
(cf.~\cite{JarPfl2000}, Th.~3.6.6).

Fix a $j\in\{1,\dots,N\}$ and a point
$$
(a',a'')\in (A_1\tdots A_{j-1})\times
(A_{j+1}\tdots A_N)\setminus\Sigma_j(Q_f).
$$
Take an
$a_j\in D_j\setminus (Q_f)_{(a',a'')}$ and let $r>0$ be such that
$\Delta_a(r)\subset\wdht X\setminus\wdht Q_f$, where $a=(a',a_j,a'')$.
Take a $D'_j\Subset D_j\setminus\wdht M_{(a',a'')}$ with $a_j\in D'_j$.
We may assume that
$\Delta_{(a',a'')}(r)\times D'_j\subset\wdht X\setminus\wdht M$ and
$\Delta_{a_j}(r)\subset D'_j$.
By the Rothstein theorem \ref{Rots} with $p:=n_1+\dots+n_{j-1}
+n_{j+1}+\dots+n_N$, $q:=n_j$,
$$
A:=((A_1\tdots A_{j-1})\times(A_{j+1}\tdots A_N))\cap\Delta_{(a',a'')}(r),
$$
we get an open set
$\Omega_a\supset A\times D'_j$ such that $\wdht f$ extends meromorphically
to $\Omega_a$.

The proof of Theorem 1.3 is completed.


\halfskip

\section{Proof of Theorem 1.4.}

\halfskip

It suffices to prove that for each $j$ there exists an open neighborhood
$\Omega_j$ of the cross $X_j:=\MX(A_j,B_j;D_j,G_j)=
(A_j\times G_j)\cup(D_j\times B_j)$ such that there exists an
$\wdtl f_j\in\CM(\Omega_j)$ with $\wdtl f_j=f$ on $X_j\setminus S$.

Indeed, we may assume that $\Omega_j\subset\wdht X_j$. Observe that
$\wdht X_j\nearrow\wdht X$. By Lemma 9 from \cite{JarPfl2002b}
the envelope of holomorphy of $\Omega_j$ equals $\wdht X_j$. Hence, by
Theorem 3.6.6 from \cite{JarPfl2000}, the function $\wdtl f_j$ extends to a function
$\wdht f_j\in\CM(\wdht X_j)$. Since $X_j\setminus S$ is not pluripolar
(by (1.4.2)),
we conclude that $\wdht f_j=\wdht f_{j+1}$ on $\wdht X_j$. Finally, we
glue up the functions $(\wdht f_j)_{j=1}^\infty$ and we get the required
extension.

\halfskip

Fix $(a,b)\in A_j\times B_j\setminus S$ and let $r>0$ be such that
$\Delta_{(a,b)}(r)\subset D_j\times G_j\setminus S$. Define
$Y:=\MX(A\cap\Delta_a(r),B\cap\Delta_b(r);\Delta_a(r),\Delta_b(r))$.
Then $f\in\CO_s(Y)$ and hence, by Theorem \ref{MT}, $f|_Y$ extends
holomorphically on $\wdht Y$. In particular, $f$ extends holomorphically
to an open neighborhood of $(a,b)$.

\halfskip

By the Rothstein theorem \ref{Rots}, we get an open set
$$
\Omega_{j,a,b}=(\Delta_a(r_{a,b})\times G_j)\cup(D_j\times
\Delta_b(r_{a,b}))\subset D_j\times G_j
$$
for which there exists
a function $\wdht f_{j,a,b}\in\CM(\Omega_{j,a,b})$ such that
$\wdht f_{j,a,b}=f$ on $X\cap\Omega_{j,a,b}\setminus S$.

\halfskip

Now we show that if $\Omega_{j,a,b}\cap\Omega_{j,a',b'}\neq\varnothing$,
then  $\wdht f_{j,a,b}=\wdht f_{j,a',b'}$ on
$\Omega_{j,a,b}\cap\Omega_{j,a',b'}$. Observe that
\begin{align*}
\Omega_{j,a,b}\cap\Omega_{j,a',b'}&=(\Delta_a(r_{a,b})\cap\Delta_{a'}(r_{a',b'}))\times G_j\\
&\cup\Delta_a(r_{a,b})\times\Delta_{b'}(r_{a',b'})\\
&\cup\Delta_{a'}(r_{a',b'})\times\Delta_b(r_{a,b})\\
&\cup D_j\times (\Delta_b(r_{a,b})\cap\Delta_{b'}(r_{a',b'})).
\end{align*}

First observe that $\wdht f_{j,a,b}=f=\wdht f_{j,a',b'}$ on
$(A_j\times B_j)\cap(\Delta_a(r_{a,b})\times\Delta_{b'}(r_{a',b'}))\setminus S$.
Hence, by (1.4.2), $\wdht f_{j,a,b}=\wdht f_{j,a',b'}$ on
$\Delta_a(r_{a,b})\times\Delta_{b'}(r_{a',b'})$. The same argument works on
$\Delta_{a'}(r_{a',b'})\times\Delta_b(r_{a,b})$.

If $\Delta_a(r_{a,b})\cap\Delta_{a'}(r_{a',b'})\neq\varnothing$, then
for any $\beta\in B_j$ we have
$\wdht f_{j,a,b}(\cdot,\beta)=f(\cdot,\beta)$ on
$A_j\cap\Delta_a(r_{a,b})\setminus S_{(\cdot,\beta)}$. Hence
$\wdht f_{j,a,b}(\cdot,\beta)=\wdtl{f(\cdot,\beta)}$ on $\Delta_a(r_{a,b})$,
and, consequently,
$\wdht f_{j,a,b}(\cdot,\beta)=\wdtl{f(\cdot,\beta)}=\wdht f_{j,a',b'}(\cdot,\beta)$
on $\Delta_a(r_{a,b})\cap\Delta_{a'}(r_{a',b'})$ for any $\beta\in B_j$. The
identity principle implies that $\wdht f_{j,a,b}=\wdht f_{j,a',b'}$
on $(\Delta_a(r_{a,b})\cap\Delta_{a'}(r_{a',b'}))\times G_j$.
The same argument works on
$D_j\times (\Delta_b(r_{a,b})\cap\Delta_{b'}(r_{a',b'}))$.

\halfskip

It remains to observe that, by (1.4.5),
$\Omega_j:=\bigcup_{(a,b)\in A_j\times B_j\setminus S}
\Omega_{j,a,b}$ is an open neighborhood of $X_j$.

The proof of Theorem 1.4 is completed.


\bibliographystyle{amsplain}

\end{document}